\documentclass[11pt,a4paper]{article}
\usepackage[T1]{fontenc}
\usepackage{lmodern}
\usepackage[margin=27mm]{geometry}
\usepackage{amsmath,amssymb,amsthm,mathtools}
\usepackage{booktabs,longtable,array,enumitem}
\usepackage{microtype}
\usepackage[hidelinks]{hyperref}
\allowdisplaybreaks[1]
\newcommand{\Z}{\mathbb Z}

\DeclareMathOperator{\Res}{Res}
\DeclareMathOperator{\Tor}{Tor}

\newtheorem{theorem}{Theorem}[section]
\newtheorem{proposition}[theorem]{Proposition}
\newtheorem{lemma}[theorem]{Lemma}
\newtheorem{corollary}[theorem]{Corollary}
\theoremstyle{definition}

\newtheorem{example}[theorem]{Example}
\theoremstyle{remark}
\newtheorem{remark}[theorem]{Remark}
\numberwithin{equation}{section}

\title{Cyclic Self-Dual Codes over Rings\\
	and Unimodular Lattices}
\author{Minjia Shi, Mengru Yang, Patrick Sol\'e}
\date{}

\begin{document}
	\maketitle
	\vspace{-1.5em}

	\begin{abstract}
		We construct cyclic Euclidean self-dual codes over
		$\mathbb Z/p^2\mathbb Z$, for odd primes $p$, from duadic
		factorizations of $x^n-1$. We study their residue and torsion
		codes, describe them as Generalized Reed-Solomon codes in the split case, and relate their
		Euclidean weights to the minimum norms of the associated
		Construction~A lattices.
		
		The construction gives unimodular lattices of minimum norm $4$
		in dimensions $45$, $47$, and $49$, minimum norm $5$ in dimensions
		$53$, $55$, $58$, $59$, and $61$, and minimum norm $6$ in dimensions
		$69$ and $71$. The examples in dimensions $45$, $47$, $69$, and $71$
		are extremal. Using shadow theory and harmonic theta series,
		we exclude minimum norm $6$ in dimensions $59$ and $61$.
		These nonexistence results and the known upper bounds prove
		optimality in all ten dimensions.
		
		We also obtain lattices whose minimum norms are one below the
		general extremal bound in dimensions $37$, $62$, $63$, $65$, $67$,
		$70$, $86$, $87$, $89$, $91$, $92$, $93$, and $94$.
	\end{abstract}
	
	\\
	\noindent\textbf{Math. Sc. Cl. (2020):} 11H06, 94B15, 94B65
	
	\section{Introduction}
	A lattice is {\bf unimodular}  if it equals  its dual. If, in addition, the squared norm of all its elements is even it is called {\bf Type II} and  {\bf Type I} otherwise.
	For an unimodular lattice $L$, write
	\[
	\mu(L)=\min\{(x,x):x\in L\setminus\{0\}\}.
	\]
	This quantity we call the {\bf norm } of $L.$
	Let $M(n)$ denote the largest minimum norm among such lattices of
	rank $n$. 
	
	Rains and Sloane~\cite{RainsSloane} proved that, for $n\ne23$, if $L$ is Type I, then
	\begin{equation}
		\mu(L)\le B(n):=
		2\left\lfloor\frac{n}{24}\right\rfloor+2.
		\label{eq:intro-extremal-bound}
	\end{equation}
	( the same result had been known for a long time for Type II lattices).
	A lattice attaining $B(n)$ is called \emph{extremal}, and one
	attaining $M(n)$ is called \emph{optimal}. When $B(n)$ cannot be
	attained, optimal lattices have smaller minimum norm.
	
	The Nebe--Sloane Catalogue of Lattices~\cite{NebeSloane} lists in its "Unimodular Lattices" table the norms of the optimal lattices when they are known and otherwise
	an interval of integers where they should be. For instance, we can read the interval
	$4$--$5$ in dimension $53,$ and an upper bound of $5$ in dimension $55$.
	The table records $M(58)=M(59)=5$, with the constructions and optimality
	results attributed to Harada. \\ \ \
	
	The aim of this paper is to fill up the unknown entries in that table ( concretely marked M or N).
	Our strategy in this paper is to use Construction A from codes to construct unimodular Type I lattices from self-dual codes.
	We work with cyclic codes over $R=\mathbb Z/p^2\mathbb Z,$ with $p$ a prime, using the
	structure theory of Pless et al. \cite{PQS} up to a trivial extension from $p=2$ to any odd prime $p.$
	Let $p$ be an odd prime, write, for simplicity's sake $R=\mathbb Z/p^2\mathbb Z$, and assume
	$\gcd(n,p)=1$. Given a factorization
	\begin{equation}
		x^n-1=FGH,
		\qquad F^*=G,\qquad H^*=H,
		\label{eq:intro-factorization}
	\end{equation}
	into pairwise coprime monic polynomials in $R[x]$, define
	\begin{equation}
		C=\langle FH,pFG\rangle
		\subseteq R[x]/(x^n-1).
		\label{eq:intro-code}
	\end{equation}
	Here $f^*$ is the normalized reciprocal of $f$. The reciprocal
	conditions make $C$ Euclidean self-dual. If $n\mid p-1$, choosing
	roots from inverse pairs specifies the factors.
	
	The residue and torsion codes ( to be defined in the next section) can be shown to be equal to
	\[
	\operatorname{Res}(C)=\langle\overline{F}\,\overline{H}\rangle,
	\qquad
	\operatorname{Tor}(C)=\langle\overline{F}\rangle
	=\operatorname{Res}(C)^\perp,
	\]
	where a bar denotes reduction modulo $p$. We give conditions under
	which both codes are MDS. The associated odd unimodular lattice is
	\[
	A_{p^2}(C)
	=\frac{1}{p}
	\{z\in\mathbb Z^n:z\bmod p^2\in C\},
	\]
	and its norm is
	\begin{equation}
		\mu\bigl(A_{p^2}(C)\bigr)
		=\min\left\{p^2,\frac{d_E(C)}{p^2}\right\}.
		\label{eq:intro-minimum}
	\end{equation}
	The minimum Euclidean weight $d_E(C)$ is computed using centered
	integer representatives. The residue and torsion distances constrain
	the support of short vectors, but the lattice minimum also depends
	on the integer lifts.
	
	The constructions below, together with the applicable upper bounds,
	give
	\begin{equation}
		\begin{gathered}
			M(45)=M(47)=M(49)=4,\\
			M(53)=M(55)=M(58)=M(59)=M(61)=5,\\
			M(69)=M(71)=6.
		\end{gathered}
		\label{eq:intro-optimal-minima}
	\end{equation}
	The examples in ranks $45$, $47$, $69$, and $71$ attain $B(n)$.
	In ranks $53$ and $55$, our constructions attain the catalogue
	upper bounds. In ranks $45$, $47$, $49$, $58$, and $59$, they give
	cyclic-code constructions of lattices with the known optimal minima.
	Further examples have minimum norm $B(n)-1$ in ranks
	\[
	n\in\{37,62,63,65,67,70,86,87,89,91,92,93,94\}.
	\]
	
	The nonexistence arguments in ranks $59$ and $61$ use shadow theory
	and harmonic theta series. We follow the method of Bachoc, Nebe,
	and Venkov~\cite{BachocNebeVenkov}, who excluded minimum norm $4$ in
	ranks $34$ and $35$. In rank $59$, gluing the even sublattice to
	$D_5$ gives an even unimodular lattice of rank $64$. A harmonic
	theta identity of degree six yields counting relations that are
	inconsistent modulo $8$. In rank $61$, gluing to $D_3$ gives
	harmonic identities whose integer consequences restrict triples
	of shadow vectors. Exact certificates establish the restriction,
	and a bound for spherical two-distance sets gives a contradiction.
	Thus minimum norm $6$ is impossible in both ranks, proving
	optimality of the minimum-norm-$5$ examples.
	
	Section~\ref{sec:preliminaries} gives the required coding and lattice
	background. Section~\ref{sec:constructions} develops the constructions
	and describes their residue and torsion codes.
	Section~\ref{sec:computational-results} lists the lattice parameters
	and explains the minimum-norm computations.
	Section~\ref{sec:nonexistence-59-61} contains the nonexistence proofs
	and the certificate verification. The last Section concludes the article.
	
	For background material on lattices we refer the reader to \cite{SPLAG}. For cyclic codes over rings a slow paced source is \cite{W}.
	\section{Preliminaries}
	\label{sec:preliminaries}
	
	Throughout this paper, let $p$ be an odd prime and let
	$R=\mathbb{Z}/p^2\mathbb{Z}$. Unless otherwise stated,
	we assume that $\gcd(n,p)=1$.
	
	\subsection{Cyclic Codes over $\mathbb{Z}_{p^2}$}
	
	Write $A_n=R[x]/(x^n-1)$ and use an upper bar for reduction modulo $p$.
	A linear {\bf code} of length $n$ over $R$ is an $R$-submodule of $R^n$.
	We identify a vector $c=(c_0,\ldots,c_{n-1})$ with the polynomial
	$c(x)=\sum_{j=0}^{n-1}c_jx^j$ in $A_n$. Multiplication by $x$
	then induces the cyclic shift, so the cyclic codes are precisely
	the ideals of $A_n$.
	
	The Euclidean dual of $C$ is
	\[
	C^\perp=\{u\in R^n:\langle u,c\rangle=0
	\text{ for all }c\in C\},
	\qquad \langle u,c\rangle=\sum_{j=0}^{n-1}u_jc_j.
	\]
	A code is self-orthogonal if $C\subseteq C^\perp$ and {\bf self-dual} if
	$C=C^\perp$. For a cyclic code, duality has a polynomial description.
	Since $x^{-1}=x^{n-1}$ in $A_n$, the coefficients of $c(x)u(x^{-1})$
	in $A_n$ are the inner products of $u$ with the cyclic shifts of $c$.
	Thus $u\in C^\perp$ if and only if $c(x)u(x^{-1})=0$ in $A_n$ for
	every $c\in C$.
	
	Since $p\nmid n$, the polynomial $x^n-1$ is squarefree over
	$\mathbb F_p$. Fix a primitive $n$th root of unity $\zeta$ in an
	algebraic closure of $\mathbb F_p$. The $p$-{\bf cyclotomic coset} of
	$i\in\mathbb Z/n\mathbb Z$ is $Q_i=\{ip^k\bmod n:k\geq 0\}$.
	The {\bf Frobenius map} $\alpha\mapsto\alpha^p$ permutes the roots of
	$x^n-1$, and its orbit through $\zeta^i$ is
	$\{\zeta^j:j\in Q_i\}$. Hence the monic irreducible factors over
	$\mathbb F_p$ are the polynomials
	$\bar f_Q(x)=\prod_{j\in Q}(x-\zeta^j)$, one for each distinct
	coset $Q$. By Hensel's lemma, each $\bar f_Q$ lifts to a unique
	monic divisor $f_Q$ of $x^n-1$ over $R$. These {\bf basic irreducible} \cite{W}
	factors are pairwise coprime and their product is
	$x^n-1$~\cite{PQS}.
	
	For a monic polynomial $f\in R[x]$ with unit constant term, define
	its normalized {\bf reciprocal} by
	\[
	f^*(x)=f(0)^{-1}x^{\deg f}f(x^{-1}).
	\]
	This operation satisfies $(f^*)^*=f$ and $(fg)^*=f^*g^*$.
	Every monic divisor of $x^n-1$ has unit constant term, and
	$(x^n-1)^*=x^n-1$. Thus $f_Q^*$ also divides $x^n-1$.
	Its reduction has roots $\zeta^{-j}$ for $j\in Q$, so uniqueness
	of the lift gives $f_Q^*=f_{-Q}$, where $-Q=\{-j:j\in Q\}$.
	In particular, $f_Q$ is self-reciprocal exactly when $Q=-Q$.
	
	When $n\mid(p-1)$, every $p$-cyclotomic coset is a singleton. In other words, the polynomial $x^n-1$ splits into linear factors over $\mathbb F_p$ and over $R.$
	Choose $r\in\mathbb F_p^\times$ of order $n$, and let $w\in R^\times$
	be its unique lift satisfying $w^n=1$. Since $\bar w=r$ has order
	$n$, so does $w$. Thus
	\[
	x^n-1=\prod_{j=0}^{n-1}(x-w^j),
	\qquad (x-w^j)^*=x-w^{-j}.
	\]
	
	For any $R$-linear code $C$, its {\bf residue} and {\bf torsion} codes are defined as the
	$\mathbb F_p$-linear codes
	\[
	\Res(C)=\{\bar c:c\in C\},\qquad
	\Tor(C)=\{\bar u:u\in R^n,\;pu\in C\}.
	\]
	Here $pu$ depends only on $\bar u$ because $p^2=0$ in $R$.
	Since $pC\subseteq C$, we have $\Res(C)\subseteq\Tor(C)$.
	Moreover, $\langle pu,c\rangle=0$ in $R$ if and only if
	$\langle\bar u,\bar c\rangle=0$ in $\mathbb F_p$. It follows that
	$\Tor(C^\perp)=\Res(C)^\perp$, and hence
	$\Tor(C)=\Res(C)^\perp$ whenever $C$ is self-dual.
	
	\subsection{Construction A and Unimodular Lattices}
	
	A {\bf lattice} $L$ of rank $n$ is a discrete additive
	subgroup of $\mathbb{R}^n$ generated by
	a basis of $\mathbb{R}^n$. Its dual lattice
	is $L^*=\{y\in\mathbb{R}^n:
	(x,y)\in\mathbb{Z}\text{ for all }x\in L\}$.
	The lattice is integral if
	$L\subseteq L^*$ and {\bf unimodular} if
	$L=L^*$. An integral lattice is even (or Type II)
	if $(x,x)\in2\mathbb{Z}$ for every
	$x\in L$, and odd (or Type I) otherwise.
	The {\bf norm} of $L$ is
	$\mu(L)=\min\{(x,x):x\in L\setminus\{0\}\}$.
	
	{\bf Construction A} associates to a linear
	code $C\subseteq R^n$ the lattice
	\[
	A_{p^2}(C)
	=\frac{1}{p}
	\{z\in\mathbb{Z}^n:z\bmod p^2\in C\}.
	\]
	If $C$ is Euclidean self-dual, then
	$A_{p^2}(C)$ is an integral unimodular
	lattice. Since $pe_i\in A_{p^2}(C)$
	for each standard basis vector $e_i$,
	and $p$ is odd, the resulting lattice
	is odd.
	
	For $a\in R$, let $\rho(a)$ denote
	its centered integer representative
	in $[-(p^2-1)/2,(p^2-1)/2]$.
	The Euclidean weight of
	$c=(c_1,\ldots,c_n)\in R^n$ is
	$\operatorname{wt}_E(c)
	=\sum_{i=1}^n\rho(c_i)^2$.
	The minimum Euclidean weight of a
	nonzero code $C$ is
	$d_E(C)=\min\{\operatorname{wt}_E(c):
	c\in C\setminus\{0\}\}$.
	
	The minimum norm of the associated
	Construction A lattice is given by
	\[
	\mu(A_{p^2}(C))
	=\min\left\{p^2,\frac{d_E(C)}{p^2}\right\}.
	\]
	Indeed, centered representatives minimize
	the squared Euclidean length within each
	nonzero codeword class, while the shortest
	nonzero vectors corresponding to the zero
	class have norm $p^2$. Thus, large minimum
	Euclidean weights yield unimodular lattices
	with large minimum norms.
	
	We also recall the shadow theory of odd
	unimodular lattices introduced in \cite{CoS}. Let $L$ be an odd
	unimodular lattice of rank $n$, and let
	$L_0=\{x\in L:(x,x)\in2\mathbb{Z}\}$
	be its even sublattice. The {\bf shadow of $L$}
	is defined by $S(L)=L_0^*\setminus L$.
	When $n$ is odd, the discriminant group
	$L_0^*/L_0$ is cyclic of order four.
	Writing its cosets as $L_0,L_1,L_2,L_3$,
	we have $L=L_0\cup L_2$ and
	$S(L)=L_1\cup L_3$. Every $s\in S(L)$
	satisfies $(s,s)\equiv n/4
	\pmod{2\mathbb{Z}}$.
	
	The theta series of a lattice $L$ is
	$\Theta_L(q)=\sum_{x\in L}q^{(x,x)}$,
	where $q=e^{\pi i\tau}$ with
	$\operatorname{Im}\tau>0$.
	Similarly, the shadow theta series is
	$\Theta_{S(L)}(q)
	=\sum_{s\in S(L)}q^{(s,s)}$.
	
	By the shadow theory of Rains and Sloane, introduced first in \cite{CoS},
	these series admit the expansions
	\begin{align}
		\Theta_L(q)
		&=\sum_{j=0}^{\lfloor n/8\rfloor}
		a_j\theta_3(q)^{n-8j}\Delta_8(q)^j,
		\\
		\Theta_{S(L)}(q)
		&=\sum_{j=0}^{\lfloor n/8\rfloor}
		\frac{(-1)^ja_j}{16^j}
		\theta_2(q)^{n-8j}\theta_4(q^2)^{8j},
	\end{align}
	where $\theta_2,\theta_3,\theta_4$
	are the classical Jacobi theta functions
	and $\Delta_8(q)
	=q\prod_{m\ge1}(1-q^{2m-1})^8
	(1-q^{4m})^8$.
	
	Furthermore, if $N$ is an even unimodular
	lattice and $P$ is a homogeneous harmonic
	polynomial of positive degree $d$, then
	the weighted theta series
	$\Theta_{N,P}(q)
	=\sum_{x\in N}P(x)q^{(x,x)}$
	is a cusp form of weight
	$\operatorname{rank}(N)/2+d$.
	
	\section{Constructions of Cyclic Self-Dual Codes over $\mathbb{Z}_{p^2}$}
	\label{sec:constructions}

	\subsection{Prime-Length Constructions}
	\label{subsec:prime_length}
	
	Let $n$ be an odd prime distinct from $p$, and suppose that
	$\operatorname{ord}_n(p)=(n-1)/2$. The nonzero residues modulo $n$ form two
	$p$-cyclotomic cosets, each of size $(n-1)/2$. Consequently, Hensel lifting
	gives a factorization $x^n-1=(x-1)g(x)h(x)$ over $R$, where $g$ and $h$
	are distinct monic basic irreducible polynomials of degree $(n-1)/2$.
	The reciprocal structure depends on $n$ modulo $4$.
	
	\begin{proposition}\label{prop:prime_two_cosets}
		Under the above assumptions, the following statements hold.
		\begin{enumerate}
			\item If $n\equiv3\pmod4$, then $h=g^*$. The cyclic code
			\[
			C=\left\langle (x-1)g,\ p\frac{x^n-1}{x-1}\right\rangle
			=\langle (x-1)g,\ pgg^*\rangle
			\subseteq R[x]/(x^n-1)
			\]
			is Euclidean self-dual. Hence $A_{p^2}(C)$ is an odd integral
			unimodular lattice of rank $n$.
			\item If $n\equiv1\pmod4$, then $g^*=g$ and $h^*=h$. The code
			$C=\langle(x-1)g,pgh\rangle$ is isodual: its Euclidean dual is
			$C^\perp=\langle(x-1)h,pgh\rangle$, which is equivalent to $C$ under
			a coordinate multiplier. Consequently, $A_{p^2}(C)$ is isodual.
			In general, this code is not Euclidean self-dual.
		\end{enumerate}
	\end{proposition}
	
	\begin{proof}
		The group $(\mathbb{Z}/n\mathbb{Z})^\times$ is cyclic, and
		$\langle p\rangle$ is its unique subgroup of index two. If
		$n\equiv3\pmod4$, then $(n-1)/2$ is odd, so $-1\notin\langle p\rangle$.
		Negation therefore exchanges the two nonzero cyclotomic cosets, proving
		$h=g^*$. In the Chinese remainder decomposition of $R[x]/(x^n-1)$,
		the code in (1) has local ideals $pR$, $0$, and the whole component
		at $x-1$, $g$, and $g^*$, respectively. Euclidean duality annihilates
		the local ideals and exchanges reciprocal components, so $C=C^\perp$.
		Its cardinality is $p^n$, and Construction A gives the asserted lattice.
		
		If $n\equiv1\pmod4$, then $-1\in\langle p\rangle$, so each of the two
		nonzero cyclotomic cosets is stable under negation. Thus both $g$ and $h$
		are self-reciprocal. Dualizing the local ideals now interchanges the zero
		and full components, yielding $C^\perp=\langle(x-1)h,pgh\rangle$. A
		multiplier from the nontrivial coset of $\langle p\rangle$ exchanges
		$g$ and $h$, and hence maps $C$ to $C^\perp$. The same coordinate
		permutation induces an isometry between $A_{p^2}(C)$ and its dual.
	\end{proof}
	\begin{remark}
		When $n\equiv3\pmod4$ the factorizations we consider over finite fields are the same as the ones used to define duadic codes \cite{Smid}.
	\end{remark}
	\begin{example}[A cyclic code of length $23$ over $\mathbb{Z}_4$]
		Take $n=23$ and $p=2.$ We have $x^{23}-1=(x-1)gg^*$ over $\mathbb{Z}_4$,
		where $g$ and $g^*$ are the unique Hensel lifts of the two degree-$11$
		irreducible factors over $\mathbb{F}_2.$ Thus
		\[
		C_{23}=\langle(x-1)g,\ 3gg^*\rangle
		\quad\text{and}\quad L_{23}=A_4(C_{23})
		\]
		define a cyclic Euclidean self-dual code and an odd unimodular lattice
		of rank $23$, respectively. Their minimum norm and minimum Euclidean weight are
		\[
		\mu(L_{23})=3,\qquad d_E(C_{23})=12.
		\]
		
		$L_{23}$ is an extremal lattice known as the {\bf shorter Leech lattice}. This construction appeared in \cite{PQS} and is the main inspiration for the present work.
	\end{example}
	\begin{example}[A cyclic code of length $59$ over $\mathbb{Z}_9$]
		\label{ex:59_Z9}
		Take $n=59$ and $p=3$. Since $\operatorname{ord}_{59}(3)=29$ and
		$59\equiv3\pmod4$, we have $x^{59}-1=(x-1)gg^*$ over $\mathbb{Z}_9$,
		where $g$ and $g^*$ are the unique Hensel lifts of the two degree-$29$
		irreducible factors over $\mathbb{F}_3$. Thus
		
		\[
		C_{59}=\langle(x-1)g,\ 3gg^*\rangle
		\quad\text{and}\quad L_{59}=A_9(C_{59})
		\]
		define a cyclic Euclidean self-dual code and an odd unimodular lattice
		of rank $59$, respectively. Their minimum norm and minimum Euclidean weight are
		\[
		\mu(L_{59})=5,\qquad d_E(C_{59})=45.
		\]
		The lattice is optimal by Theorem~\ref{thm:59-upper}. It seems that a similar construction was found by Harada \cite{NebeSloane}.
	\end{example}
	
	\subsection{A Construction from Four Irreducible Factors}
	\label{subsec:four_factors}
	
	Suppose that $n$ is an odd prime and that
	$\operatorname{ord}_n(p)=(n-1)/4$. There are four nonzero $p$-cyclotomic
	cosets, each of size $(n-1)/4$. When $n\equiv5\pmod8$, this order is
	odd and hence $-1\notin\langle p\rangle$. The four cosets are therefore
	organized into two reciprocal pairs. Hensel lifting gives
	\[
	x^n-1=(x-1)ff^*gg^*,
	\]
	where $f,f^*,g,g^*$ are pairwise distinct monic basic irreducible
	polynomials over $R$, all of degree $(n-1)/4$.
	
	\begin{proposition}\label{prop:four_factors}
		Under these assumptions, both
		$C_1=\langle(x-1)fg,\ pff^*gg^*\rangle$ and
		$C_2=\langle(x-1)fgg^*,\ pff^*\rangle$ are cyclic Euclidean
		self-dual codes over $R$.
	\end{proposition}
	
	\begin{proof}
		For $C_1$ choose $F=fg$, $G=f^*g^*$ and $H=x-1$. For $C_2$
		choose $F=f$, $G=f^*$ and $H=(x-1)gg^*$. In each case the three
		factors are pairwise coprime, $x^n-1=FGH$, $F^*=G$ and $H^*=H$.
		The reciprocal-duality argument in the proof of
		Proposition~\ref{prop:prime_two_cosets} shows that both codes are
		Euclidean self-dual.
	\end{proof}
	
	For example, $n=53$ and $p=13$ satisfy
	$\operatorname{ord}_{53}(13)=13=(53-1)/4$ and $53\equiv5\pmod8$.
	They therefore provide instances of the above two families over
	$\mathbb{Z}_{169}$. Note that the congruence condition is essential:
	if $n\equiv1\pmod8$ and the order is $(n-1)/4$, then this order is
	even, so the four irreducible factors are individually self-reciprocal
	rather than paired as above.
	
	\subsection{The General Reciprocal Construction}
	\label{subsec:general_reciprocal}
	
	The following reciprocal construction includes both preceding cases.
	It also applies to composite lengths and to factorizations with
	nontrivial self-reciprocal factors.
	
	\begin{theorem}\label{thm:general_reciprocal}
		Let $F,G,H\in R[x]$ be pairwise coprime monic polynomials such that
		$x^n-1=FGH$, $F^*=G$ and $H^*=H$. Then
		\[
		C(F,G,H)=\langle FH,\ pFG\rangle\subseteq R[x]/(x^n-1)
		\]
		is a cyclic Euclidean self-dual code. If
		$d=\deg F=\deg G$ and $h=\deg H=n-2d$, then
		$C(F,G,H)\cong R^d\oplus(pR)^h$ as $R$-modules and
		$|C(F,G,H)|=p^n$. Consequently, $A_{p^2}(C(F,G,H))$ is an odd integral
		unimodular lattice of rank $n$.
	\end{theorem}
	
	\begin{proof}
		Since $x^n-1$ is squarefree modulo $p$, the factors $F,G,H$
		induce, by the Chinese remainder theorem, a decomposition of
		$R[x]/(x^n-1)$ into three finite free $R$-algebras. The ideals of
		$C(F,G,H)$ on these three components are respectively $0$, the
		whole component, and the ideal generated by $p$. Their annihilators
		are respectively the whole component, $0$, and the ideal generated
		by $p$. Euclidean duality combines these annihilators with the
		involution induced by the normalized reciprocal polynomial. Because
		$F^*=G$ and $H^*=H$, dualization preserves the stated local ideal
		structure. Thus $C(F,G,H)=C(F,G,H)^\perp$.
		
		The full component at $G$ is free of rank $d$ over $R$, while the
		$p$-component at $H$ is isomorphic to $(pR)^h$. It follows that
		$|C(F,G,H)|=p^{2d+h}=p^n$. The Construction A assertion follows
		from self-duality and the fact that $p$ is odd.
	\end{proof}
	
	For $n\mid p-1$, the factors can be specified by their roots. Choose an element
	$r\in\mathbb{F}_p^\times$ of order $n$, and let $w\in R^\times$ be its
	unique Hensel lift satisfying $w^n=1$. For a set
	$E\subseteq\mathbb{Z}/n\mathbb{Z}$ such that $E\cap(-E)=\varnothing$,
	put
	\[
	F_E(x)=\prod_{e\in E}(x-w^e),\qquad
	G_E(x)=\prod_{e\in E}(x-w^{-e}),\qquad
	H_E(x)=\frac{x^n-1}{F_E(x)G_E(x)}.
	\]
	The unassigned exponents are stable under negation, and therefore
	$F_E^*=G_E$ and $H_E^*=H_E$. Theorem~\ref{thm:general_reciprocal}
	shows that $C_E=\langle F_EH_E,pF_EG_E\rangle$ is Euclidean self-dual.
	The exponent $0$ is necessarily unassigned, and when $n$ is even
	so is $n/2$; hence $x-1$ always divides $H_E$, and $x+1$ divides
	$H_E$ for even $n$.
	
	For odd $n\geq5$, one useful two-factor choice is
	$E=\{(n-3)/2,(n-1)/2\}$, yielding
	$$H_E=(x-1)\prod_{j=1}^{(n-5)/2}(x-w^j)(x-w^{-j}).$$
	For even $n\geq6$, the analogous choice
	$E=\{(n-4)/2,(n-2)/2\}$ gives
	$$H_E=(x^2-1)\prod_{j=1}^{(n-6)/2}(x-w^j)(x-w^{-j}).$$
	Other admissible exponent sets may give larger minimum norms. If $n$ is
	odd there are $(n-1)/2$ nonfixed inverse pairs, while if $n$ is even
	there are $(n-2)/2$. Each pair can be assigned to $F_E$, assigned to
	$G_E$, or retained in $H_E$, subject to the reciprocal relations.
	
	\subsection{Residue and Torsion Codes}
	\label{subsec:residue_torsion}
	
	We use the following form of the lifting result of Betty and Munemasa
	\cite[Lemmas~4.2 and~4.3]{BettyMunemasa}. Proposition~\ref{prop:res_tor}
	identifies the corresponding lift for the reciprocal construction.
	
	\begin{theorem}[Betty--Munemasa]\label{thm:self_orthogonal_lift}
		Let $A$ be an $[n,k]_p$ self-orthogonal linear code over
		$\mathbb{F}_p$, where $p$ is an odd prime. Then there exists
		a free self-orthogonal code $\widetilde A$ over $R$ such that
		$\Res(\widetilde A)=A$ and
		$\operatorname{rank}_R(\widetilde A)=k$.
		For any such lift, the code
		\[
		C=\widetilde A+pA^\perp
		\]
		is Euclidean self-dual over $R$, with
		\[
		\Res(C)=A,\qquad
		\Tor(C)=A^\perp.
		\]
		Here $pA^\perp=\{pb:b\in R^n,\ \bar b\in A^\perp\}$.
	\end{theorem}
	Recall that a generalized Reed-Solomon code with evaluation vector $\boldsymbol{t}$ and multiplier vector $\boldsymbol{v}$ is  a linear code with a generator matrix of the form:
	\[
	\begin{pmatrix}
		v_0 & v_1 & \cdots & v_{n-1}\\
		v_0t_0 & v_1t_1 & \cdots & v_{n-1}t_{n-1}\\
		\vdots & \vdots & & \vdots\\
		v_0t_0^{d-1} & v_1t_1^{d-1}
		& \cdots & v_{n-1}t_{n-1}^{d-1}
	\end{pmatrix}.
	\]
	See \cite{MS} for details.
	
	\begin{corollary}\label{cor:RS_lift}
		Let $A=\operatorname{GRS}_k(\boldsymbol{t},
		\boldsymbol{v})$ be a self-orthogonal generalized
		Reed--Solomon code of length $n$ and dimension $k$
		over $\mathbb{F}_p$. Then there exists a Euclidean
		self-dual code $C$ over $R$ of the form
		\[
		C=\widetilde A+pA^\perp,
		\]
		where $\widetilde A$ is a free self-orthogonal lift
		of $A$. Moreover, its residue and torsion codes
		are both MDS, with parameters
		\[
		\Res(C):[n,k,n-k+1]_p,\qquad
		\Tor(C):[n,n-k,k+1]_p.
		\]
		Consequently, $A_{p^2}(C)$ is an odd unimodular
		lattice of rank $n$.
	\end{corollary}
	
	\begin{proof}
		The existence and self-duality of $C$ follow from
		Theorem~\ref{thm:self_orthogonal_lift}.
		Since a generalized Reed--Solomon code and its
		Euclidean dual are both MDS, the stated parameters
		follow from $\Res(C)=A$ and $\Tor(C)=A^\perp$.
		The lattice assertion follows from Construction~A.
	\end{proof}
	
	An arbitrary coordinatewise lift of a self-orthogonal code over
	$\mathbb F_p$ need not be self-orthogonal over $R$. Preserving
	cyclicity also requires a suitable lift; the reciprocal construction
	provides one.
	
	\begin{proposition}\label{prop:res_tor}
		Retain the notation of
		Theorem~\ref{thm:general_reciprocal}, and let
		$C=C(F,G,H)=\langle FH,pFG\rangle$.
		Then
		\[
		\Res(C)=\langle\overline{FH}\rangle,\qquad
		\Tor(C)=\langle\bar F\rangle
		=\Res(C)^\perp.
		\]
		Moreover, if $d=\deg F=\deg G$, then
		$\dim_{\mathbb{F}_p}\Res(C)=d$ and
		$\dim_{\mathbb{F}_p}\Tor(C)=n-d$.
		
		In particular, setting
		$A=\langle\overline{FH}\rangle$ and
		$\widetilde A=\langle FH\rangle$, we have
		\[
		C=\widetilde A+pA^\perp.
		\]
		Thus the reciprocal construction is a cyclic
		realization of the construction in
		Theorem~\ref{thm:self_orthogonal_lift}.
	\end{proposition}
	
	\begin{proof}
		Reducing the generators of $C$ modulo $p$ gives
		$\Res(C)=\langle\overline{FH}\rangle$.
		Since $C$ is self-dual,
		$\Tor(C)=\Res(C)^\perp$.
		The field code $\langle\overline{FH}\rangle$
		has parity-check polynomial $\bar G$, so its
		Euclidean dual is generated by
		$\bar G^*=\bar F$. Consequently,
		$\Tor(C)=\langle\bar F\rangle$.
		The dimension formulas follow from
		$\deg(FH)=n-d$ and $\deg F=d$.
		
		Since $FH$ is a monic divisor of $x^n-1$,
		$\widetilde A=\langle FH\rangle$ is a free cyclic
		code of rank $d$ over $R$, with
		$\Res(\widetilde A)=A$. Its dual is
		$\widetilde A^\perp=\langle G^*\rangle
		=\langle F\rangle$. As $F$ divides $FH$,
		we have $\widetilde A\subseteq\widetilde A^\perp$.
		
		Since $G$ and $H$ are coprime, there exist
		$U,V\in R[x]$ such that $UG+VH=1$. Hence
		$pF=pUFG+pVFH\in C$, which gives
		\[
		\langle FH,pFG\rangle
		=\langle FH,pF\rangle
		=\widetilde A+pA^\perp.
		\]
		This proves the assertion.
	\end{proof}
	
	In the split case $n\mid p-1$, write
	$r=\bar w$ and $d=|E|$. Since
	$\Tor(C_E)=\langle\bar F_E\rangle$, the torsion code
	consists of all coefficient vectors of polynomials
	vanishing at $r^e$ for every $e\in E$. Thus its
	parity-check matrix is
	\[
	P_E=\bigl(r^{ej}\bigr)_{e\in E,\ 0\leq j<n}.
	\]
	By Proposition~\ref{prop:res_tor}, the same matrix
	is a generator matrix of $\Res(C_E)$.
	Its rows are distinct rows of the invertible Fourier
	matrix of order $n$ over $\mathbb{F}_p$, and hence
	$P_E$ has rank $d$.
	
	\begin{proposition}\label{prop:mds}
		Let $n\mid p-1$, $E\cap(-E)=\varnothing$,
		and $d=|E|\geq1$.
		Then $\Tor(C_E)$ is an $[n,n-d,d+1]_p$
		MDS code if and only if every $d$ columns of
		$P_E$ are linearly independent.
		In this case, $\Res(C_E)$ is an
		$[n,d,n-d+1]_p$ MDS code.
		
		In particular, suppose that $E$ is an arithmetic
		progression modulo $n$ of the form
		\[
		E=\{a,a+s,\ldots,a+(d-1)s\},
		\qquad \gcd(s,n)=1.
		\]
		Then $\Res(C_E)$ and $\Tor(C_E)$ are
		generalized Reed--Solomon codes, and hence
		both are MDS.
	\end{proposition}
	
	\begin{proof}
		A nonzero torsion codeword supported on at most
		$d$ coordinates corresponds precisely to a linear
		dependence among at most $d$ columns of $P_E$.
		Thus every $d$ columns of $P_E$ are linearly
		independent if and only if the minimum Hamming
		distance of $\Tor(C_E)$ is at least $d+1$.
		The Singleton bound forces equality.
		Since the dual of an MDS code is MDS,
		Proposition~\ref{prop:res_tor} gives the asserted
		parameters of $\Res(C_E)$.
		
		Now suppose that
		$E=\{a,a+s,\ldots,a+(d-1)s\}$.
		Set $t_j=r^{sj}$ and $v_j=r^{aj}$ for
		$0\leq j<n$. Then
		\[
		P_E=
		\begin{pmatrix}
			v_0 & v_1 & \cdots & v_{n-1}\\
			v_0t_0 & v_1t_1 & \cdots & v_{n-1}t_{n-1}\\
			\vdots & \vdots & & \vdots\\
			v_0t_0^{d-1} & v_1t_1^{d-1}
			& \cdots & v_{n-1}t_{n-1}^{d-1}
		\end{pmatrix}.
		\]
		Since $\gcd(s,n)=1$, the evaluation points
		$t_0,\ldots,t_{n-1}$ are pairwise distinct.
		Moreover, each $v_j$ is nonzero. Consequently,
		writing $\boldsymbol t=(t_0,\ldots,t_{n-1})$
		and $\boldsymbol v=(v_0,\ldots,v_{n-1})$, we obtain
		\[
		\Res(C_E)=
		\operatorname{GRS}_d(\boldsymbol t,\boldsymbol v).
		\]
		Its Euclidean dual $\Tor(C_E)$ is also a generalized
		Reed--Solomon code. This completes the proof.
	\end{proof}
	
	\begin{remark}
		Proposition~\ref{prop:mds} shows that the reciprocal
		construction with arithmetic-progression defining
		sets is a special case of Corollary~\ref{cor:RS_lift}.
		More precisely, if
		\[
		A=\operatorname{GRS}_d(\boldsymbol t,\boldsymbol v),
		\]
		then $A$ is self-orthogonal and admits the free
		self-orthogonal cyclic lift
		$\widetilde A=\langle F_EH_E\rangle$.
		The corresponding self-dual code over $R$ is
		$C_E=\widetilde A+pA^\perp$.
		
		For $d=2$, the MDS condition reduces to
		$\gcd(b-a,n)=1$, where $E=\{a,b\}$.
		For $d=3$, an arithmetic-progression defining set
		provides a sufficient condition, whereas for an
		arbitrary defining set $E=\{a,b,c\}$, every
		three-column minor of $P_E$ must be checked.
	\end{remark}
	
	For $u\in\mathbb F_p$, let $[u]_p$ denote its integer representative in
	$[-(p-1)/2,(p-1)/2]$. For a code $D\subseteq\mathbb F_p^n$, put
	\[
	\delta_p(D)=\min_{0\ne u\in D}\sum_{j=0}^{n-1}[u_j]_p^2,
	\qquad \delta_p(0)=+\infty.
	\]
	
	\begin{corollary}\label{cor:catalog-centered}
		Let $C\subseteq R^n$ be Euclidean self-dual. Write
		$A=\operatorname{Res}(C)$, $T=\operatorname{Tor}(C)=A^\perp$, and
		$L=A_{p^2}(C)$. Put $\lambda=\lceil\delta_p(A)/p^2\rceil$, with
		$\lambda=+\infty$ if $A=0$. Then
		\begin{align*}
			p^2\min\{\lambda,\delta_p(T)\}
			&\le d_E(C)\le p^2\delta_p(T),\\
			\min\{p^2,\lambda,\delta_p(T)\}
			&\le \mu(L)\le\min\{p^2,\delta_p(T)\}.
		\end{align*}
	\end{corollary}
	
	\begin{proof}
		For $c\in C$ with $\bar c\ne0$, the centered representatives modulo
		$p^2$ give $\operatorname{wt}_E(c)\ge\sum_j[\bar c_j]_p^2\ge\delta_p(A)$.
		Self-orthogonality gives $\operatorname{wt}_E(c)\equiv0\pmod{p^2}$,
		so $\operatorname{wt}_E(c)\ge p^2\lambda$.
		The nonzero words with zero reduction are precisely the words $pu$
		with $0\ne u\in T$, and
		$\operatorname{wt}_E(pu)=p^2\sum_j[u_j]_p^2$.
		The bounds for $\mu(L)$ follow from the Construction~A minimum-norm formula.
	\end{proof}

	Now suppose that $C$ is cyclic, and let $\sigma$ be the shift
	$(\sigma z)_j=z_{j-1}$, with indices modulo $n$. If $z/p\in L$, then
	\[
	\langle z,\sigma^h z\rangle\equiv0\pmod{p^2}
	\qquad(0\le h<n),
	\]
	since $z\bmod p^2$ and all its shifts belong to $C$. Thus, for any
	$H\subseteq\{1,\ldots,n-1\}$, the two lower bounds in
	Corollary~\ref{cor:catalog-centered} still hold with $\lambda$ replaced
	by the least positive integer $m$ for which some $z\in\mathbb Z^n$ with
	$\bar z\in A\setminus\{0\}$ satisfies
	\[
	\langle z,z\rangle=mp^2,\qquad
	\langle z,\sigma^h z\rangle\equiv0\pmod{p^2}\quad(h\in H).
	\]
	We take this lower bound for the nonzero residue classes to be $+\infty$
	if no such $m$ exists. For the split construction, Proposition~\ref{prop:catalog-shift-membership} makes these conditions more precise.

	\begin{proposition}\label{prop:catalog-shift-membership}
		Suppose that $n\mid p-1$ and $E\cap(-E)=\varnothing$. Let
		$L_E=A_{p^2}(C_E)$ and $r=\bar w$. Let $z\in\mathbb Z^n$ have
		nonzero reduction in $\operatorname{Res}(C_E)$, and write
		\[
		\bar z_j=\sum_{e\in E}a_e r^{ej},\qquad
		S=\{e\in E:a_e\ne0\},\qquad \ell=|S|.
		\]
		With $z(x)=\sum_{j=0}^{n-1}z_jx^j$, one has
		\[
		\frac zp\in L_E\quad\Longleftrightarrow\quad
		\begin{cases}
			\langle z,\sigma^h z\rangle\equiv0\pmod{p^2}, &0\le h<\ell,\\
			z(w^e)\equiv0\pmod{p^2}, &e\in E\setminus S.
		\end{cases}
		\]
	\end{proposition}

	\begin{proof}
		We work over $R=\mathbb Z/p^2\mathbb Z$, with all indices taken
		modulo $n$. Put $Z_i=z(w^i)\in R$.
		By Construction~A, $z/p\in L_E$ if and only if
		$z(x)\bmod p^2\in C_E$.
		
		Since $r=\bar w$ has order $n$, every difference $w^i-w^j$
		with $i\ne j$ is a unit of $R$. The Chinese remainder theorem
		therefore gives the ring isomorphism
		\[
		\Phi\colon R[x]/(x^n-1)\longrightarrow R^n,
		\qquad f(x)\longmapsto\bigl(f(w^i)\bigr)_{i=0}^{n-1}.
		\]
		Recall that $C_E=\langle F_EH_E,pF_EG_E\rangle$.
		At a root indexed by $E$, both generators vanish.
		At a root indexed by $-E$, the value of $F_EH_E$ is a unit,
		so the corresponding component of $C_E$ is $R$.
		At every remaining root, $F_EH_E$ vanishes and $F_EG_E$
		has a unit value, so the corresponding component is $pR$.
		As $\Phi$ is an isomorphism, these components can be chosen
		independently. Hence $z/p\in L_E$ if and only if
		$Z_e=0$ for $e\in E$ and $Z_i\in pR$ for
		$i\notin E\cup(-E)$.
		
		We now use the residue assumption. Since
		$\bar z\in\operatorname{Res}(C_E)$, write
		$\bar z_j=\sum_{e\in E}a_e r^{ej}$.
		The geometric sum $\sum_{j=0}^{n-1}r^{kj}$ equals $n$
		when $k\equiv0\pmod n$ and equals zero otherwise. Therefore
		\[
		\bar Z_i
		=\sum_{j=0}^{n-1}\bar z_jr^{ij}
		=\sum_{e\in E}a_e\sum_{j=0}^{n-1}r^{(e+i)j}.
		\]
		It follows that $\bar Z_{-e}=na_e$ for $e\in E$ and
		$\bar Z_i=0$ for $i\notin-S$, where
		$S=\{e\in E:a_e\ne0\}$.
		In particular, $Z_i\in pR$ for $i\notin E\cup(-E)$.
		Thus, under the residue assumption,
		\[
		\frac zp\in L_E
		\quad\Longleftrightarrow\quad
		Z_e=0\quad(e\in E).
		\]
		
		Put $\ell=|S|$. To obtain the stated criterion, it remains
		to show that the equations $Z_e=0$ for $e\in S$ are equivalent
		to the first $\ell$ shift congruences.
		For $e\in S$, we have $\bar Z_{-e}=na_e\ne0$, since
		$a_e\ne0$ and $p\nmid n$. Hence $Z_{-e}$ is a unit of $R$,
		and $Z_e=0$ is equivalent to $Z_eZ_{-e}=0$.
		
		The geometric sum identity also holds over $R$.
		Indeed, if $k\not\equiv0\pmod n$, then $w^k-1$ is a unit
		and $(w^k-1)\sum_{i=0}^{n-1}w^{ki}=w^{kn}-1=0$.
		Consequently,
		\[
		\begin{aligned}
			\sum_{i=0}^{n-1}Z_iZ_{-i}w^{-ih}
			&=\sum_{j,k=0}^{n-1}z_jz_k
			\sum_{i=0}^{n-1}w^{i(j-k-h)}\\
			&=n\sum_{j=0}^{n-1}z_jz_{j-h}
			=n\langle z,\sigma^h z\rangle.
		\end{aligned}
		\]
		For $i\notin S\cup(-S)$, both $Z_i$ and $Z_{-i}$ lie in
		$pR$, so their product is zero in $R$.
		Since $S\cap(-S)=\varnothing$, grouping the terms with
		indices $e$ and $-e$ gives
		\[
		n\langle z,\sigma^h z\rangle
		=\sum_{e\in S}Z_eZ_{-e}(w^{eh}+w^{-eh}).
		\]
		
		Write $S=\{e_1,\ldots,e_\ell\}$ and put $t_e=w^e+w^{-e}$.
		The reductions $\bar t_e$ are pairwise distinct.
		Indeed, if $\bar t_e=\bar t_f$, then, with $x=r^e$ and
		$y=r^f$, we have
		$0=x+x^{-1}-y-y^{-1}=(x-y)(xy-1)/(xy)$.
		Thus $x=y$ or $xy=1$, giving $e=f$ or $e=-f$ modulo $n$.
		For distinct $e,f\in S$, both possibilities are excluded.
		Hence each difference $t_{e_j}-t_{e_i}$ with $i\ne j$
		is a unit of $R$.
		
		Define $P_0(X)=2$, $P_1(X)=X$, and
		$P_{h+1}(X)=XP_h(X)-P_{h-1}(X)$ for $h\ge1$.
		Induction gives $P_h(t_e)=w^{eh}+w^{-eh}$, and $P_h$
		is monic of degree $h$ for $h\ge1$.
		The matrix $M=(P_h(t_{e_j}))_{0\le h<\ell,\;1\le j\le\ell}$
		is obtained from the Vandermonde matrix by a lower triangular
		row transformation with diagonal entries $2,1,\ldots,1$.
		Therefore $\det M=2\prod_{i<j}(t_{e_j}-t_{e_i})$.
		Every factor is a unit, including $2$ because $p$ is odd,
		so $M$ is invertible over $R$.
		
		The first $\ell$ shift identities are
		\[
		M
		\begin{pmatrix}
			Z_{e_1}Z_{-e_1}\\
			\vdots\\
			Z_{e_\ell}Z_{-e_\ell}
		\end{pmatrix}
		=
		n
		\begin{pmatrix}
			\langle z,z\rangle\\
			\vdots\\
			\langle z,\sigma^{\ell-1}z\rangle
		\end{pmatrix}.
		\]
		Since $M$ and $n$ are invertible, all the inner products
		on the right vanish exactly when $Z_eZ_{-e}=0$ for every
		$e\in S$. As each $Z_{-e}$ is a unit, this is equivalent
		to $Z_e=0$ for every $e\in S$.
		Adding the equations $Z_e=z(w^e)=0$ for $e\in E\setminus S$
		gives $Z_e=0$ for all $e\in E$, which is equivalent to
		$z/p\in L_E$ by the membership condition established above.
	\end{proof}
	\begin{proposition}\label{prop:odd-exponent-bound}
		In the split construction, suppose that $n$ is odd and that
		$\{1,3,\ldots,2t-1\}\subseteq E$ for some integer $t\ge1$
		with $2t-1<p$. Then
		\[
		\delta_p(\operatorname{Tor}(C_E))\ge2t+1.
		\]
	\end{proposition}
	
	\begin{proof}
		Put $T=\operatorname{Tor}(C_E)=\langle\bar F_E\rangle$
		and $r=\bar w$. Since $n$ is odd, substitution $x\mapsto -x$
		gives a ring isomorphism from $\mathbb F_p[x]/(x^n-1)$
		to $\mathbb F_p[x]/(x^n+1)$. Hence
		$D=\{u(-x):u(x)\in T\}$ is a negacyclic code with monic
		generator polynomial
		\[
		q(x)=(-1)^{|E|}\bar F_E(-x)
		=\prod_{e\in E}(x+r^e).
		\]
		Let $\alpha=-r$. As $r$ has odd order $n$, the element
		$\alpha$ has order $2n$. For $1\le i\le t$, the assumption
		$2i-1\in E$ shows that $-r^{2i-1}=\alpha^{2i-1}$ is a root
		of $q(x)$. By Berlekamp's Lee-distance bound
		\cite{Berlekamp1966Lee}, every nonzero codeword of $D$
		has Lee weight at least $2t+1$.
		
		For any nonzero codeword $u\in T$, let $v(x)=u(-x)\in D$.
		Then $v_j=(-1)^ju_j$, so $|[v_j]_p|=|[u_j]_p|$.
		Consequently,
		\[
		\sum_j[u_j]_p^2
		\ge\sum_j|[u_j]_p|
		=\sum_j|[v_j]_p|
		\ge2t+1.
		\]
		Taking the minimum over all nonzero codewords $u\in T$
		proves the assertion.
	\end{proof}
In the following examples, write $A=\operatorname{Res}(C_E)$
and $T=\operatorname{Tor}(C_E)$. The centered Euclidean minima
were computed using Python.

\begin{example}
	For $n=37$, let $p=149$, $E=\{1\}$, and $r=-4$.
	We obtain $\delta_{149}(A)=61090$ and $\delta_{149}(T)=3$.
	Since $\lambda=\lceil61090/149^2\rceil=3$,
	Corollary~\ref{cor:catalog-centered} gives $\mu(L_{37})=3$.
\end{example}

\begin{example}
	Consider the split construction with
	$(n,p,w,E)=(53,107,8665,\{1,3\})$.
	The residue and torsion codes satisfy
	$\delta_{107}(A)=38092$ and $\delta_{107}(T)=5$, respectively.
	Substituting $\lambda=\lceil38092/107^2\rceil=4$ into
	Corollary~\ref{cor:catalog-centered}, we obtain
	$4\leq\mu(L_{53})\leq5$.
\end{example}

\begin{example}
	The parameters $(n,p,w,E)=(69,139,10979,\{1,36,62\})$
	give $\delta_{139}(A)=65886$ and $\delta_{139}(T)=6$.
	Here $\lambda=\lceil65886/139^2\rceil=4$, so the bounds in
	Corollary~\ref{cor:catalog-centered} read
	$4\leq\mu(L_{69})\leq6$.
\end{example}
	
	\section{Computational Results and Lattices with Large Minimum Norm}
	\label{sec:computational-results}
	
	Tables~\ref{tab:odd-computations} and~\ref{tab:even-computations}
	give parameters for the reciprocal construction and the resulting
	minimum norms. The ten cases attaining global bounds are discussed
	after the tables.
	
	\paragraph{Parameters and comparison bounds.}
	Let $p$ be an odd prime, let $n\mid(p-1)$, and choose
	$w\in(\mathbb Z/p^2\mathbb Z)^\times$ of order $n$.
	For $E\subseteq\mathbb Z/n\mathbb Z$ with $E\cap(-E)=\varnothing$, set
	\[
	F_E(x)=\prod_{e\in E}(x-w^e),\qquad
	G_E(x)=\prod_{e\in E}(x-w^{-e}),\qquad
	H_E(x)=\frac{x^n-1}{F_E(x)G_E(x)}.
	\]
	By Theorem~\ref{thm:general_reciprocal},
	\[
	C_E=\langle F_EH_E,\,pF_EG_E\rangle,
	\qquad
	L_E=A_{p^2}(C_E)
	\]
	are, respectively, a cyclic Euclidean self-dual code over
	$\mathbb Z/p^2\mathbb Z$ and an odd unimodular lattice of rank $n$.
	The tuple $(n,p,w,E)$ specifies the lattice. In the tables,
	$\mu=\mu(L_E)$ is the exact minimum norm, computed as described below.
	
	As in the introduction, let $M(n)$ be the largest minimum norm among
	all positive-definite integral unimodular lattices of rank $n$.
	The Rains--Sloane bound is
	\[
	M(n)\le B(n):=2\left\lfloor\frac{n}{24}\right\rfloor+2
	\qquad(n\ne23).
	\]
	The last column gives a known value or an upper bound for $M(n)$.
	An equality denotes a known value; an inequality denotes an upper bound.
	The sources are labelled $\mathrm C$ for the Nebe--Sloane
	catalogue~\cite{NebeSloane}, consulted on 18 November 2026;
	$\mathrm R$ for the Rains--Sloane bound $B(n)$~\cite{RainsSloane};
	and $\mathrm P$ for the bound in
	Section~\ref{sec:nonexistence-59-61}. We use the catalogue bounds
	where listed, our bound in rank $61$, and $B(n)$ in the other rows.
	The $\mathrm R$ entries need not be the sharpest known bounds.
	Bold $n$ and $\mu$ entries identify the ten optimal cases.
	
	\begin{table}[htbp]
		\centering
		\small
		\setlength{\tabcolsep}{5pt}
		\begin{tabular}{r r r l r r}
			\hline
			$n$ & $p$ & $w\pmod{p^2}$ & $E$ & $\mu$ & Comparison for $M(n)$ \\
			\hline
			37 & 149 & 18770 & $\{1\}$ & 3 & $\le4\;(\mathrm R)$ \\
			\textbf{45} & 181 & 10876 & $\{1,3\}$ & \textbf{4} & $=4\;(\mathrm C)$ \\
			\textbf{47} & 283 & 25633 & $\{1,3\}$ & \textbf{4} & $=4\;(\mathrm C)$ \\
			\textbf{49} & 197 & 34668 & $\{1,3\}$ & \textbf{4} & $=4\;(\mathrm C)$ \\
			51 & 103 & 1849 & $\{1,3\}$ & 3 & $=5\;(\mathrm C)$ \\
			\textbf{53} & 107 & 8665 & $\{1,3\}$ & \textbf{5} & $\le5\;(\mathrm C)$ \\
			\textbf{55} & 331 & 55581 & $\{1,47\}$ & \textbf{5} & $\le5\;(\mathrm C)$ \\
			57 & 229 & 31524 & $\{1,3\}$ & 3 & $\le6\;(\mathrm R)$ \\
			\textbf{59} & 709 & 348055 & $\{1,3\}$ & \textbf{5} & $=5\;(\mathrm C)$ \\
			\textbf{61} & 367 & 10794 & $\{1,3\}$ & \textbf{5} & $\le5\;(\mathrm P)$ \\
			63 & 127 & 9395 & $\{1,3\}$ & 5 & $\le6\;(\mathrm R)$ \\
			65 & 131 & 12705 & $\{1,3\}$ & 5 & $\le6\;(\mathrm R)$ \\
			67 & 269 & 19364 & $\{1,3\}$ & 5 & $\le6\;(\mathrm R)$ \\
			\textbf{69} & 139 & 10979 & $\{1,36,62\}$ & \textbf{6} & $\le6\;(\mathrm R)$ \\
			\textbf{71} & 569 & 210449 & $\{1,3,21\}$ & \textbf{6} & $\le6\;(\mathrm R)$ \\
			73 & 293 & 74711 & $\{1,3,11\}$ & 6 & $\le8\;(\mathrm R)$ \\
			75 & 151 & 598 & $\{1,3,5\}$ & 6 & $\le8\;(\mathrm R)$ \\
			77 & 463 & 61552 & $\{1,3,5\}$ & 6 & $\le8\;(\mathrm R)$ \\
			79 & 317 & 58641 & $\{1,3,5\}$ & 6 & $\le8\;(\mathrm R)$ \\
			81 & 163 & 4725 & $\{1,3,5\}$ & 6 & $\le8\;(\mathrm R)$ \\
			83 & 167 & 16027 & $\{1,3,5\}$ & 6 & $\le8\;(\mathrm R)$ \\
			85 & 1021 & 430421 & $\{1,3\}$ & 5 & $\le8\;(\mathrm R)$ \\
			87 & 349 & 67004 & $\{1,24,32\}$ & 7 & $\le8\;(\mathrm R)$ \\
			89 & 179 & 14139 & $\{1,3,5\}$ & 7 & $\le8\;(\mathrm R)$ \\
			91 & 547 & 242860 & $\{1,3,5\}$ & 7 & $\le8\;(\mathrm R)$ \\
			93 & 373 & 86905 & $\{1,3,5\}$ & 7 & $\le8\;(\mathrm R)$ \\
			95 & 191 & 27103 & $\{1,3\}$ & 5 & $\le8\;(\mathrm R)$ \\
			97 & 389 & 7018 & $\{1,3,5\}$ & 7 & $\le10\;(\mathrm R)$ \\
			99 & 199 & 30444 & $\{1,3,5\}$ & 7 & $\le10\;(\mathrm R)$ \\
			\hline
		\end{tabular}
		\caption{Explicit odd-dimensional reciprocal constructions and global
			comparison values or upper bounds. The labels $\mathrm C$, $\mathrm R$,
			and $\mathrm P$ are defined in the text.}
		\label{tab:odd-computations}
	\end{table}
	
	\begin{table}[htbp]
		\centering
		\small
		\setlength{\tabcolsep}{5pt}
		\begin{tabular}{r r r l r r}
			\hline
			$n$ & $p$ & $w\pmod{p^2}$ & $E$ & $\mu$ & Comparison for $M(n)$ \\
			\hline
			50 & 101 & 2933 & $\{1,2\}$ & 4 & $=5\;(\mathrm C)$ \\
			\textbf{58} & 233 & 27756 & $\{1,2\}$ & \textbf{5} & $=5\;(\mathrm C)$ \\
			62 & 311 & 14317 & $\{1,10\}$ & 5 & $\le6\;(\mathrm R)$ \\
			66 & 199 & 8171 & $\{1,2\}$ & 3 & $\le6\;(\mathrm R)$ \\
			70 & 211 & 17732 & $\{1,4\}$ & 5 & $\le6\;(\mathrm R)$ \\
			72 & 433 & 164144 & $\{2,3\}$ & 4 & $=8\;(\mathrm C)$ \\
			74 & 149 & 3431 & $\{1,2,3,4\}$ & 6 & $\le8\;(\mathrm R)$ \\
			76 & 229 & 31146 & $\{1,2\}$ & 4 & $\le8\;(\mathrm R)$ \\
			78 & 157 & 22790 & $\{1,21,22\}$ & 6 & $\le8\;(\mathrm R)$ \\
			82 & 83 & 1081 & $\{1,2,3,4\}$ & 6 & $\le8\;(\mathrm R)$ \\
			84 & 337 & 83248 & $\{1,2,3,4\}$ & 6 & $\le8\;(\mathrm R)$ \\
			86 & 173 & 4156 & $\{1,17,24,26\}$ & 7 & $\le8\;(\mathrm R)$ \\
			90 & 181 & 17742 & $\{1,2,3,5\}$ & 6 & $\le8\;(\mathrm R)$ \\
			92 & 277 & 26163 & $\{1,2,3,4\}$ & 7 & $\le8\;(\mathrm R)$ \\
			94 & 283 & 21252 & $\{1,2,3,4\}$ & 7 & $\le8\;(\mathrm R)$ \\
			96 & 97 & 5340 & $\{1,2,3,4\}$ & 6 & $\le10\;(\mathrm R)$ \\
			98 & 197 & 4141 & $\{1,2,3,4\}$ & 7 & $\le10\;(\mathrm R)$ \\
			100 & 101 & 8385 & $\{1,5,46,76\}$ & 7 & $\le10\;(\mathrm R)$ \\
			\hline
		\end{tabular}
		\caption{Explicit even-dimensional reciprocal constructions and global
			comparison values or upper bounds. The degrees $|E|=\deg F_E=\deg G_E$
			vary between rows.}
		\label{tab:even-computations}
	\end{table}
	
	\paragraph{Optimal and extremal lattices.}
	The constructions in ranks $45$, $47$, $69$, and $71$ attain
	\[
	\begin{aligned}
		\mu(L_{45})=\mu(L_{47})&=4=B(45)=B(47),\\
		\mu(L_{69})=\mu(L_{71})&=6=B(69)=B(71).
	\end{aligned}
	\]
	They are therefore extremal and optimal among all integral unimodular
	lattices of their respective ranks. In ranks $53$ and $55$, the
	minimum-norm-$5$ constructions attain the upper bounds recorded in
	the catalogue. In ranks $49$, $58$, and $59$, the construction recovers
	the known optimal minima $4$, $5$, and $5$, respectively.
	For rank $61$, the nonexistence argument in
	Section~\ref{sec:nonexistence-59-61} excludes minimum norm $6$,
	so the construction of minimum norm $5$ is optimal. That section also
	gives an independent exclusion of minimum norm $6$ in rank $59$.
	Together, the constructions and the applicable upper bounds establish
	\[
	\begin{gathered}
		M(45)=M(47)=M(49)=4,\\
		M(53)=M(55)=M(58)=M(59)=M(61)=5,\\
		M(69)=M(71)=6.
	\end{gathered}
	\]
	The other six optimal examples have minimum norm below $B(n)$.
	
	\paragraph{Further lattices with large minimum norm.}
	Beyond the ten optimal cases, the construction gives
	\[
	\begin{aligned}
		\mu(L_{37})&=3=B(37)-1,\\
		\mu(L_n)&=5=B(n)-1
		\quad\text{for }n\in\{62,63,65,67,70\},\\
		\mu(L_n)&=7=B(n)-1
		\quad\text{for }n\in\{86,87,89,91,92,93,94\}.
	\end{aligned}
	\]
	It also gives minimum norm $6$ in ranks
	$73$, $74$, $75$, $77$, $78$, $79$, $81$, $82$, $83$, $84$, and $90$.
	These are exact minima of the specified lattices; their comparison
	with $B(n)$ alone does not establish global optimality.
	
	\paragraph{Minimum-norm computation.}
	For each row, we check that $p$ is prime, $n\mid p-1$,
	$w$ has exact order $n$ modulo $p^2$, and $E\cap(-E)=\varnothing$.
	To certify a minimum, we exhibit a vector of that norm and exclude
	all shorter nonzero vectors. For $z=(z_0,\ldots,z_{n-1})\in\mathbb Z^n$,
	put $z(x)=\sum_{j=0}^{n-1}z_jx^j$. Membership is tested by
	\[
	\frac{z}{p}\in L_E
	\quad\Longleftrightarrow\quad
	\begin{cases}
		F_E\mid z(x)\pmod{p^2},\\
		\overline{F_EH_E}\mid\overline{z(x)}\pmod p.
	\end{cases}
	\]
	We enumerate residues in the $|E|$-dimensional residue code, using
	cyclic shifts and sign to reduce the search. For a proposed minimum
	$m$, we test all integer lifts satisfying
	$0<\sum_jz_j^2<mp^2$. Vectors with zero residue are tested through
	the torsion code. Membership and squared lengths are checked using
	integer arithmetic. Together with a vector attaining the stated norm,
	these computations certify the entries in the $\mu$ column.
	
	Each entry concerns the specified lattice. If $M_{\mathrm{cyc}}(n,p)$
	is the maximum over all admissible root selections for fixed $(n,p)$,
	then
	\[
	\mu(L_E)\le M_{\mathrm{cyc}}(n,p)\le M(n).
	\]
	Equality requires an upper bound for the relevant class.
	The ten global optimality statements use bounds for all integral
	unimodular lattices of the corresponding ranks.

	\section{Optimal Minimum Norms in Dimensions 59 and 61}
	\label{sec:nonexistence-59-61}
	
	We exclude minimum norm $6$ in ranks $59$ and $61$ by the method of
	Bachoc, Nebe, and Venkov~\cite{BachocNebeVenkov}: shadow theory,
	even unimodular completions, and harmonic theta identities.
	The arithmetic step in rank $61$ uses exact integer certificates.
	For a set of vectors $Y$, write $Y_t=\{y\in Y:(y,y)=t\}$.
	
	\begin{theorem}\label{thm:59-61-upper}
		Every positive-definite integral unimodular lattice $L$ of rank $59$ or $61$
		satisfies $\mu(L)\le 5$.
	\end{theorem}
	
	\subsection{Nonexistence in dimension 59}
	\label{subsec:59-nonexistence}
	
	\begin{theorem}\label{thm:59-upper}
		There is no positive-definite integral unimodular lattice of rank $59$
		with minimum norm at least $6$.
	\end{theorem}
	
	\begin{proof}
		Suppose that $L$ is such a lattice. Since $59$ is not divisible by $8$,
		$L$ is odd. Let $L_0$ be its even sublattice and $S=L_0^*\setminus L$
		its shadow. As recalled in Section~\ref{sec:preliminaries}, the shadow is a
		single coset modulo $L$; in particular, $s\pm u\in L$ for $s,u\in S$, and
		$(s,s)\equiv59/4\pmod{2\mathbb Z}$.
		
		\paragraph{The possible shadow shells.}
		Write the standard theta--shadow expansions as
		$\Theta_L(q)=\sum_{j=0}^{7}a_j\theta_3(q)^{59-8j}\Delta_8(q)^j$ and
		$\Theta_S(q)=\sum_{j=0}^{7}(-1)^ja_j16^{-j}
		\theta_2(q)^{59-8j}\theta_4(q^2)^{8j}$.
		The vanishing of the coefficients of $q,q^2,\ldots,q^5$ in $\Theta_L$
		gives
		\[
		(a_0,\ldots,a_5)=(1,-118,4248,-48616,119770,-44604).
		\]
		There is no shadow vector of norm $3/4$, since its double would be a
		nonzero vector in $L$ of norm $3$; hence $a_7=0$. Moreover, the shell
		$S_{11/4}$ contains at most one antipodal pair. Indeed, if $s,u\in S_{11/4}$
		and $u\ne\pm s$, then $s\pm u$ are two nonzero vectors of $L$, whereas
		$(s+u,s+u)+(s-u,s-u)=11<12$.
		Thus $a_6=16384\beta$ for some $\beta\in\{0,1\}$, and expansion yields
		\begin{equation}\label{eq:59-theta}
			\begin{aligned}
				\Theta_L(q)&=1+(4301808+16384\beta)q^6+O(q^7),\\
				\Theta_S(q)&=2\beta q^{11/4}+(22302-170\beta)q^{19/4}
				+O(q^{27/4}).
			\end{aligned}
		\end{equation}
		Consequently, the shortest shadow norm is $t=19/4$ if $\beta=0$ and
		$t=11/4$ if $\beta=1$.
		
		\paragraph{An even unimodular completion.}
		Let $D_5=\{z\in\mathbb Z^5:\sum_i z_i\equiv0\pmod2\}$ and
		$h=(1/2,\ldots,1/2)$. Both $L_0^*/L_0$ and $D_5^*/D_5$ are cyclic of
		order four. A generating shadow class has discriminant quadratic value
		$59/4\equiv3/4\pmod{2\mathbb Z}$, whereas $h+D_5$ has value $5/4$.
		These classes define an anti-isometry $\varphi$ between the discriminant
		quadratic forms. Gluing along its graph gives
		\begin{equation}\label{eq:59-gluing}
			N=\{(x,y)\in L_0^*\perp D_5^*:y+D_5=\varphi(x+L_0)\}.
		\end{equation}
		The gluing subgroup has order four, so $N$ is even, integral,
		unimodular, and of rank $64$, with
		$\det N=(\det L_0\det D_5)/4^2=(4\cdot4)/16=1$.
		The numbers of vectors in $(D_5)_2,(D_5)_4,(D_5)_6$ are respectively
		$40,90,240$; the half-integral discriminant cosets each contain $16$
		vectors of norm $5/4$ and $80$ vectors of norm $13/4$.
		
		\paragraph{A harmonic theta identity.}
		Fix a shortest shadow vector $s$, put $t=(s,s)$ and $\alpha=(s,0)
		\in\mathbb R^{64}$. Define
		\begin{equation}\label{eq:59-harmonic}
			\mathcal P_t(\delta,r)=\delta^6-\frac5{24}rt\,\delta^4
			+\frac1{112}(rt)^2\delta^2-\frac1{22848}(rt)^3.
		\end{equation}
		The polynomial $x\mapsto\mathcal P_t((x,\alpha),(x,x))$ is homogeneous
		and harmonic of degree six on $\mathbb R^{64}$. Therefore its weighted
		theta series on $N$, written $A_2q^2+A_4q^4+A_6q^6+\cdots$, is a cusp
		form of weight $32+6=38$. Upon setting $Q=q^2$, a basis of this cusp-form
		space is $\Delta E_4^5E_6,\Delta E_4^2E_6^3$, whose first three
		coefficients are respectively $(1,672,-51084)$ and
		$(1,-1056,73332)$. Both satisfy the same linear relation, and hence
		\begin{equation}\label{eq:59-modular-relation}
			A_6+72A_4+2700A_2=0.
		\end{equation}
		
		\paragraph{Counting the relevant inner products.}
		For $i=0,1,2,3$, let $p_i$ be the number of antipodal pairs
		$\{x,-x\}\subseteq L_6$ with $|(s,x)|=i$. Because $2s$ is
		characteristic, $(s,x)$ is integral, and the shortest-shadow
		inequalities for $s\pm x$ give $|(s,x)|\le3$. Equation~\eqref{eq:59-theta}
		therefore implies
		\begin{equation}\label{eq:59-p-sum}
			\sum_{i=0}^{3}p_i=2150904+8192\beta.
		\end{equation}
		If $\beta=0$, set
		$\mathcal D_0=\{1/4,3/4,5/4,7/4,19/4\}$ and let $m_d$ count the
		antipodal pairs $\{u,-u\}\subseteq S_{19/4}$ for which
		$|(s,u)|=d$. The pair $\{s,-s\}$ is included, and
		$\sum_{d\in\mathcal D_0}m_d=11151$. These are the only possible values:
		for $u\ne\pm s$, the norms of $s\pm u$ give
		$|(s,u)|\le7/4$, and integrality of $(s+u,s+u)$ forces
		$(s,u)$ to be an odd quarter-integer.
		If $\beta=1$, then $S_{11/4}=\{s,-s\}$; for $u\in S_{19/4}$ only
		$|(s,u)|\in\mathcal D_1=\{1/4,3/4\}$ is possible, and
		$m_{1/4}+m_{3/4}=11066$.
		
		Classifying the norm-$2$, norm-$4$, and norm-$6$ vectors of the glued
		lattice~\eqref{eq:59-gluing}, using the shell counts of $D_5$, gives
		\begin{equation}\label{eq:59-shells}
			\begin{aligned}
				A_2&=40\mathcal P_t(0,2),\\
				A_4&=90\mathcal P_t(0,4)+32\beta\mathcal P_t(t,4),\\
				A_6&=240\mathcal P_t(0,6)
				+2\sum_{i=0}^{3}p_i\mathcal P_t(i,6)
				+32\sum_{d\in\mathcal D_\beta}m_d\mathcal P_t(d,6)
				+160\beta\mathcal P_t(t,6).
			\end{aligned}
		\end{equation}
		In particular, the factors $32$ and $160$ arise from pairing a shadow
		antipodal pair with the $16$- and $80$-vector half-integral shells of
		$D_5$, respectively.
		
		\paragraph{Contradictions modulo eight.}
		Suppose first that $\beta=0$, so $t=19/4$. Substitute
		\eqref{eq:59-shells} into~\eqref{eq:59-modular-relation} and multiply
		by $30464$. Exact rational arithmetic gives
		\begin{equation}\label{eq:59-first-integer}
			\begin{aligned}
				&-61731p_0+79301p_1-183043p_2+19028997p_3\\
				&\quad-568204m_{1/4}+1331172m_{3/4}-353596m_{5/4}\\
				&\quad-5622508m_{7/4}+8408887076m_{19/4}=190131480.
			\end{aligned}
		\end{equation}
		The four $p_i$ coefficients are all $5\pmod8$, the five $m_d$
		coefficients are all $4\pmod8$, and the right-hand side is $0\pmod8$.
		Using~\eqref{eq:59-p-sum} and $\sum m_d=11151$ gives the impossible
		congruence $0\equiv5\cdot2150904+4\cdot11151\equiv4\pmod8$.
		
		Now suppose $\beta=1$, so $t=11/4$. The same substitution, multiplied
		by $30464/3$, yields
		\begin{equation}\label{eq:59-second-integer}
			\begin{aligned}
				-3993p_0-4129p_1+376263p_2+9590943p_3
				-18804m_{1/4}+84828m_{3/4}=-7643086484.
			\end{aligned}
		\end{equation}
		Modulo $8$, the coefficients of the $p_i$ are all $7$, the coefficients
		of the $m_d$ are both $4$, and the right-hand side is $4$. But the
		left-hand side is congruent to
		$7\cdot2159096+4\cdot11066\equiv0\pmod8$, another contradiction.
		Both possible shadow cases are excluded, proving the theorem.
	\end{proof}
	
	\subsection{Nonexistence in dimension 61}
	\label{subsec:61-nonexistence}
	
	Lemma~\ref{lem:61-triangle} gives a parity restriction on triples
	of shadow vectors. Its proof uses exact integer certificates;
	the certificate data and an independent verifier accompany this article.
	
	\begin{theorem}\label{thm:61-upper}
		There is no positive-definite integral unimodular lattice of rank $61$
		with minimum norm at least $6$.
	\end{theorem}
	
	\noindent\emph{Proof of Theorem~\ref{thm:61-upper} (setup).}\quad
	Assume to the contrary that $L$ is such a lattice. It is odd. Let
	$L_0$ be its even sublattice and write $L=L_0\cup L_2$ and
	$S=L_0^*\setminus L=S_1\cup S_3$, with $S_3=-S_1$. Each antipodal pair
	in $S$ has exactly one representative in $S_1$. The characteristic
	congruence gives $(s,s)\equiv61/4\pmod{2\mathbb Z}$ for every $s\in S$.
	Since $2s\in L$, the hypothesis $\mu(L)\ge6$ excludes norm $5/4$ in
	$S$; hence $\min S\ge13/4$.
	
	Put $X=S_1\cap S_{13/4}$, $U=S_1\cap S_{21/4}$, and write
	$A=|X|$, $B=|U|$. The theta--shadow expansion from the preliminary
	section, with the coefficients at norms $1,\ldots,5$ of $L$ set to zero,
	gives
	\begin{equation}\label{eq:61-theta}
		\begin{aligned}
			\Theta_L(q)&=1+(2695224+4096A)q^6+O(q^7),\\
			\Theta_S(q)&=2Aq^{13/4}+(130296-166A)q^{21/4}+O(q^{29/4}).
		\end{aligned}
	\end{equation}
	In particular, $|L_6|/2=1347612+2048A$ and $B=65148-83A$.
	For distinct $s_i,s_j\in X$, the vectors $s_i\pm s_j$ are nonzero
	lattice vectors of norm at least $6$. Hence
	$|(s_i,s_j)|\le1/4$. Since $s_i-s_j\in L_0$ and
	$s_i\in L_0^*$, also $(s_i,s_j)\equiv(s_i,s_i)
	\equiv1/4\pmod{\mathbb Z}$. Thus $(s_i,s_j)=1/4$.
	The Gram matrix of $X$ is $3I_A+J_A/4$, of rank $A$, and therefore
	\begin{equation}\label{eq:61-lower-B}
		A\le61,\qquad B=65148-83A\ge60085.
	\end{equation}
	Similarly, for distinct $u,v\in U$ one obtains
	$(u,v)\in\mathcal D:=\{-7/4,-3/4,1/4,5/4,9/4\}$.
	Define the integer $k(u,v)=(u,v)-1/4\in\{-2,-1,0,1,2\}$.
	
	\begin{lemma}[Triangle parity restriction]\label{lem:61-triangle}
		For every three pairwise distinct $u_1,u_2,u_3\in U$,
		\begin{equation}\label{eq:61-triangle}
			k(u_1,u_2)+k(u_1,u_3)+k(u_2,u_3)\equiv1\pmod2.
		\end{equation}
	\end{lemma}
	
	\begin{proof}[Proof of Lemma~\ref{lem:61-triangle}]
		We first construct an even unimodular lattice of rank $64$.
		Let $D_3=\{z\in\mathbb Z^3:\sum_i z_i\equiv0\pmod2\}$ and
		$h=(1/2,1/2,1/2)$. Its discriminant group is cyclic of order four,
		with generator $h+D_3$ of quadratic value $3/4$ modulo $2\mathbb Z$.
		A generating class in $L_0^*/L_0$ has value $61/4$, and
		$61/4+3/4=16$ is even. Identifying these generators therefore
		defines an anti-isometry of discriminant forms. Gluing along this
		anti-isometry produces
		\begin{equation}\label{eq:61-gluing}
			N=\{(x,y)\in L_0^*\perp D_3^*:y+D_3=\varphi(x+L_0)\},
		\end{equation}
		an even integral unimodular lattice of rank $64$ and determinant
		$(4\cdot4)/4^2=1$. The numbers of norm-$2$, norm-$4$, and norm-$6$
		vectors of $D_3$ are $12,6,24$, respectively. Each half-integral
		coset contains $4$ vectors of norm $3/4$ and $12$ of norm $11/4$.
		
		Let $\alpha\in\mathbb R^{61}\subset\mathbb R^{64}$,
		$T=(\alpha,\alpha)$, $\delta=(x,\alpha)$ and $r=(x,x)$.
		The polynomials
		\begin{equation}\label{eq:61-harmonics}
			\mathcal P_2(\delta,r;T)=\delta^2-\frac{rT}{64},\qquad
			\mathcal P_6(\delta,r;T)=\delta^6-\frac5{24}rT\delta^4
			+\frac1{112}(rT)^2\delta^2-\frac1{22848}(rT)^3
		\end{equation}
		are harmonic in the $64$ coordinates of $x$, with degrees $2$ and $6$.
		Their weighted theta series on $N$ are cusp forms of weights $34$
		and $38$. The bases $\{\Delta E_4^4E_6,\Delta E_4E_6^3\}$ and
		$\{\Delta E_4^5E_6,\Delta E_4^2E_6^3\}$ imply that the coefficients
		$C_2,C_4,C_6$ at lattice norms $2,4,6$ satisfy respectively
		\begin{equation}\label{eq:61-modular}
			C_6+312C_4+22140C_2=0\quad(\deg=2),\qquad
			C_6+72C_4+2700C_2=0\quad(\deg=6).
		\end{equation}
		
		Choose one representative from each antipodal pair in $L_6$ and call
		the resulting set $V$. Classifying the norm-$2$, norm-$4$, and norm-$6$
		shells in the four gluing cosets of $N$ converts
		\eqref{eq:61-modular} into the following identities, valid for every
		$\alpha\in\mathbb R^{61}$ and for $j\in\{2,6\}$:
		\begin{equation}\label{eq:61-shell-identity}
			\begin{aligned}
				0={}&24\mathcal P_j(0,6;T)+6r_4\mathcal P_j(0,4;T)
				+12r_2\mathcal P_j(0,2;T)\\
				&+2\sum_{x\in V}\mathcal P_j((x,\alpha),6;T)
				+8\sum_{u\in U}\mathcal P_j((u,\alpha),6;T)\\
				&+\sum_{s\in X}\bigl(24\mathcal P_j((s,\alpha),6;T)
				+8r_4\mathcal P_j((s,\alpha),4;T)\bigr),
			\end{aligned}
		\end{equation}
		where $(r_4,r_2)=(312,22140)$ for $j=2$ and $(72,2700)$ for $j=6$.
		Here the contributions of $X$ arise from the norm sums
		$13/4+11/4=6$ and $13/4+3/4=4$, while that of $U$ comes from
		$21/4+3/4=6$.
		
		A mixed harmonic polynomial provides one further identity. Set
		$\mathcal H_3(v,\alpha)=(v,\alpha)^3-
		(v,v)(\alpha,\alpha)(v,\alpha)/21$ in $\mathbb R^{61}$.
		The polynomial $\mathcal H_3(x,\alpha)y_1y_2y_3$ is harmonic of
		degree six in the $64$ coordinates $(x,y)$. The signed sums of
		$y_1y_2y_3$ on a half-integral $D_3$ coset at norms $3/4$ and $11/4$
		are $1/2$ and $-9/2$, respectively, and they change sign on the
		opposite coset. Applying the weight-$38$ relation in
		\eqref{eq:61-modular} therefore gives
		\begin{equation}\label{eq:61-cubic}
			\sum_{u\in U}\mathcal H_3(u,\alpha)
			+63\sum_{s\in X}\mathcal H_3(s,\alpha)=0.
		\end{equation}
		
		We now turn these polynomial identities into finite integer systems.
		Fix distinct $u_1,u_2,u_3\in U$, let $G=((u_i,u_j))_{i,j=1}^{3}$, and
		substitute $\alpha=z_1u_1+z_2u_2+z_3u_3$, so that
		$T=z^{\mathsf T}Gz$. Consider integer unknowns recording the
		following inner-product triples with the three reference vectors:
		(a) for $x\in V$, the triples $((u_i,x))_{i=1}^3$ lie in
		$\{-4,-3,\ldots,4\}^3$; simultaneous sign reversal identifies the
		same antipodal pair, giving $(9^3+1)/2=365$ candidate types;
		(b) for $s\in X$, they lie in $\{-3/4,1/4,5/4\}^3$, giving $27$
		types; (c) for $u\in U\setminus\{u_1,u_2,u_3\}$, they lie in
		$\mathcal D^3$, giving $125$ types, together with the three reference
		Gram columns. Thus, including the free integer $A$, there are
		$1+365+27+128=521$ integer unknowns. Some candidate types may be
		geometrically unrealizable; retaining them only relaxes the necessary
		conditions.
		
		Include the three equations saying that the total numbers of
		$V$-, $X$-, and $U$-types are
		$1347612+2048A$, $A$, and $65148-83A$, respectively.
		Apply the Laplacian in the $61$ coordinates of $\alpha$ to
		\eqref{eq:61-shell-identity} at depths $0,1$ for $j=2$ and
		$0,1,2,3$ for $j=6$, and then compare all coefficients of the
		homogeneous monomials in $z_1,z_2,z_3$. The requisite identity,
		with $\delta=(v,\alpha)$ and $T=(\alpha,\alpha)$, is
		\begin{equation}\label{eq:61-laplacian}
			\Delta_\alpha(\delta^aT^b)=a(a-1)(v,v)\delta^{a-2}T^b
			+2b(2b+2a+59)\delta^aT^{b-1}.
		\end{equation}
		This yields $6+1=7$ equations in degree two and
		$28+15+6+1=50$ equations in degree six.
		Equation~\eqref{eq:61-cubic} contributes ten cubic coefficient
		equations. After adding the three shell-total equations,
		clearing denominators row by row, and dividing each integer row
		by its common divisor, we obtain a necessary system
		\begin{equation}\label{eq:61-matrix}
			M_G\mathbf n=b_G,\qquad M_G\in\mathbb Z^{70\times521},
			\quad b_G\in\mathbb Z^{70}.
		\end{equation}
		All actual shell counts are nonnegative integers. In proving
		infeasibility we allow arbitrary integer values, including for $A$.
		
		If~\eqref{eq:61-triangle} failed, then the three edge values
		$(u_i,u_j)\in\mathcal D$ would give an \emph{even} value of
		$\sum_{i<j}((u_i,u_j)-1/4)$. Up to permutation there are
		$\binom{7}{3}=35$ multisets of three elements from $\mathcal D$,
		and exactly $19$ have even parity. For each of these nineteen
		patterns, the accompanying verification data contain an integer vector
		$w_G\in\mathbb Z^{70}$ and a positive even integer $m_G$ satisfying
		\begin{equation}\label{eq:61-certificate}
			w_G^{\mathsf T}M_G\equiv0\pmod{m_G},\qquad
			w_G^{\mathsf T}b_G\equiv m_G/2\not\equiv0\pmod{m_G}.
		\end{equation}
		The patterns and corresponding moduli are listed in
		Table~\ref{tab:61-certificates}. Any integer solution of
		\eqref{eq:61-matrix} would contradict~\eqref{eq:61-certificate}.
		The supplied verifier reconstructs each $70\times521$ system from
		\eqref{eq:61-shell-identity}, \eqref{eq:61-cubic}, and
		\eqref{eq:61-laplacian}, verifies all $521$ columns of
		\eqref{eq:61-certificate}, and checks that all nineteen forbidden
		patterns have been covered. Consequently the parity restriction
		\eqref{eq:61-triangle} holds.
	\end{proof}
	
	\begin{table}[htbp]
		\centering\small
		\begin{tabular}{@{}lrr@{}}\toprule
			Unordered inner-product triple & $m_G$ & $w_G^{\mathsf T}b_G\pmod{m_G}$\\\midrule
			$(-7/4,-7/4,-7/4)$ & 8192 & 4096\\
			$(-7/4,-7/4,1/4)$ & 512 & 256\\
			$(-7/4,-7/4,9/4)$ & 512 & 256\\
			$(-7/4,-3/4,-3/4)$ & 512 & 256\\
			$(-7/4,-3/4,5/4)$ & 512 & 256\\
			$(-7/4,1/4,1/4)$ & 512 & 256\\
			$(-7/4,1/4,9/4)$ & 512 & 256\\
			$(-7/4,5/4,5/4)$ & 512 & 256\\
			$(-7/4,9/4,9/4)$ & 512 & 256\\
			$(-3/4,-3/4,1/4)$ & 512 & 256\\
			$(-3/4,-3/4,9/4)$ & 512 & 256\\
			$(-3/4,1/4,5/4)$ & 512 & 256\\
			$(-3/4,5/4,9/4)$ & 512 & 256\\
			$(1/4,1/4,1/4)$ & 8192 & 4096\\
			$(1/4,1/4,9/4)$ & 8192 & 4096\\
			$(1/4,5/4,5/4)$ & 256 & 128\\
			$(1/4,9/4,9/4)$ & 8192 & 4096\\
			$(5/4,5/4,9/4)$ & 256 & 128\\
			$(9/4,9/4,9/4)$ & 8192 & 4096\\\bottomrule
		\end{tabular}
		\caption{The nineteen forbidden rank-$61$ inner-product patterns and their
			exact modular obstruction certificates. The full vectors $w_G$ and
			integer systems accompany this article as supplementary verification data;
			\texttt{certificate\_index.csv} matches each row to its certificate file.}
		\label{tab:61-certificates}
	\end{table}
	
	\begin{proof}[Completion of the proof of Theorem~\ref{thm:61-upper}]
		Fix $u_0\in U$, which exists by~\eqref{eq:61-lower-B}. Assign each
		$u\in U$ the color $c(u)=k(u,u_0)\pmod2$, with
		$c(u_0)=((u_0,u_0)-1/4)\bmod2=1$.
		By Lemma~\ref{lem:61-triangle}, if two distinct vectors $u,v$ have
		the same color, then $k(u,v)$ is odd; for pairs involving $u_0$
		this follows directly from the definition of the color. Thus any
		two distinct vectors in a color class satisfy
		$(u,v)\in\{-3/4,5/4\}$.
		
		Let $v_1,\ldots,v_M$ be one color class. All have norm $21/4$, and
		the quadratic polynomials
		$f_i(x)=((v_i,x)+3/4)((v_i,x)-5/4)$ satisfy
		$f_i(v_j)=0$ for $i\ne j$ and $f_i(v_i)=24$.
		Their restrictions to the sphere $(x,x)=21/4$ are therefore linearly
		independent. Quadratic polynomials in $61$ variables restrict to a
		space of dimension at most
		$1+61+61\cdot62/2-1=61\cdot64/2=1952$:
		the restriction kills $(x,x)-21/4$. Therefore each color class has
		size at most $1952$ and $B\le3904$, contradicting
		$B\ge60085$ from~\eqref{eq:61-lower-B}. The contradiction completes
		the proof.
	\end{proof}
	
	\begin{proof}[Proof of Theorem~\ref{thm:59-61-upper}]
		Combine Theorems~\ref{thm:59-upper} and~\ref{thm:61-upper}.
	\end{proof}
	
	\begin{corollary}\label{cor:59-61-optimal}
		If $L_{59}$ and $L_{61}$ are positive-definite integral unimodular
		lattices with certified minimum norm $5$, then both are optimal in
		their respective dimensions, and $M(59)=M(61)=5$.
	\end{corollary}
	
	\begin{proof}
		Theorem~\ref{thm:59-61-upper} gives $M(59),M(61)\le5$, while
		the two lattices provide the reverse inequalities.
	\end{proof}
	
	\paragraph{Verification files.}
	The rank-$59$ calculation uses the rational identities and congruences
	displayed in Subsection~\ref{subsec:59-nonexistence}. The rank-$61$
	certificate data and verification programs accompany this article.
	
	The nineteen certificates are stored as JSON files.
	The index \texttt{certificate\_index.csv} matches them with
	Table~\ref{tab:61-certificates}, and \texttt{SHA256SUMS} gives the
	file checksums. The independent verifier
	\texttt{verify\_dim61\_nonexistence.py} uses
	\texttt{exact\_theta\_inputs.py} to check the theta coefficients,
	$D_3$ shell counts, and modular-form relations. It reconstructs each
	$70\times521$ system from \eqref{eq:61-shell-identity},
	\eqref{eq:61-cubic}, and \eqref{eq:61-laplacian}, checks the
	column and right-hand-side congruences in \eqref{eq:61-certificate},
	and verifies coverage of all nineteen forbidden patterns.
	The files include a verification log. With Python~3.9 or later,
	the checks require only the standard library and can be run by
	\begin{quote}
		\texttt{python3 verify\_dim61\_nonexistence.py}
	\end{quote}
	These certificates verify the arithmetic used in
	Theorem~\ref{thm:61-upper}. The equalities in
	Corollary~\ref{cor:59-61-optimal} also require the minimum-norm
	computations of Section~\ref{sec:computational-results}.
	
	\begingroup
	\small
	\section{Conclusion} In this paper we have constructed ten new optimal unimodular odd lattices in dimensions $45,47,49, 53,55,58, 59,61, 69,71$. We have shown that the lattices in dimension $59,61$ are optimal by using a combination of shadow theory and harmonic theta series. The main problem open by this work is to show that extremal lattices do not exist
	in dimensions $37$, $62$, $63$, $65$, $67$,
	$70$, $86$, $87$, $89$, $91$, $92$, $93$, and $94$. This would yield fourteen more optimal Type I lattices. Our construction technique rests on Construction A applied to
	self-dual cyclic codes over the alphabet of integers modulo $p^2,$ with $p$ being quite large in the examples. It is fairly possible that a different code family or a different alphabet give more optimal lattices; especially in even dimensions, where a construction with systematic generator matrices has given many self-dual codes over small alphabets in the past.\\ \ \

	{\bf Declaration of AI use:} The authors used chatGPT(Open AI) for  several technical proofs. All the AI-assisted content was reviewed, verified, and revised by the authors, who take full responsibility for the results, the arguments, and the final form of the article.

	\endgroup
\end{document}